\documentclass[12pt]{amsart}
\usepackage{amsmath, amsfonts, amssymb}
\newcommand{\dal}{\square}
\newcommand{\R}{{\mathbf R}}
\newcommand{\ve}{\varepsilon}
\newcommand{\pa}{\partial}

\newcommand{\jb}[1]{\left\langle #1 \right\rangle}

\newcommand{\norm}[2]{\|#1 \!:\! #2\|}
\DeclareMathOperator{\supp}{\rm supp}
\newtheorem{theorem}{Theorem}[section]

\newtheorem{lemma}[theorem]{Lemma}

\newtheorem{definition}[theorem]{Definition}

\renewcommand{\thelemma}{\thesection.\arabic{lemma}}

\numberwithin{equation}{section}

\title[Global existence for nonlinear wave equations]
{An elementary proof of global existence for nonlinear wave equations
\\
in an exterior domain}

\author[S.~Katayama]{Soichiro Katayama}
\address{Department of Mathematics, Wakayama University, 930 Sakaedani, Wakayama 640-8510, Japan}
\email{katayama@center.wakayama-u.ac.jp}
\author[H.~Kubo]{Hideo Kubo}
\address{
Department of Mathematics, Graduate School of Science,
\\
Osaka University, Toyonaka, Osaka 560-0043, Japan
}
\email{kubo@math.sci.osaka-u.ac.jp}
\thanks{The first and the second author were
partially supported by Grant-in-Aid for Young Scientists (B)
(No.~16740094), MEXT, and by  Grant-in-Aid for Science Research (No.17540157), JSPS, respectively}

\date{}

\subjclass{Primary 35L70; Secondary 35L20}

\begin{document}

\begin{abstract}
The aim of this article is to present an elementary proof of a global existence
result for nonlinear wave equations 
in an exterior domain.
The novelty of our proof is to avoid completely the scaling operator
which would make the argument complicated in the mixed problem, by using 
new weighted pointwise
estimates of a tangential derivative to the light cone.
\end{abstract}

\maketitle

\section{Introduction}
Let $\Omega$ be an unbounded domain in $\R^3$
with compact and smooth boundary $\partial \Omega$.
We put ${\mathcal O}:=\R^3 \setminus {\Omega}$, which is called
an obstacle. 
This paper is concerned with the mixed problem for 
a system of nonlinear wave equations in $\Omega$\,:
\begin{align}\label{ap1}
& (\partial_t^2-c_i^2\Delta) u_i =F_i(u, \partial u,\nabla_{\!x}\,\partial u),
& (t,x) \in (0,\infty)\times \Omega,
\\ \label{ap2}
& u(t,x)=0, & (t,x) \in (0,\infty)\times \partial\Omega,
\\ \label{ap3}
& u(0,x)=\varepsilon \phi(x),\ (\partial_t u)(0,x)=\varepsilon \psi(x), 
&  x\in \Omega,
\end{align}
for $i=1, \dots, N$, where
$c_i$ ($1\le i\le N$) are given positive constants,
$u=(u_1, \dots, u_N)$, 
$\varepsilon$ is a positive parameter and $\phi$, 
$\psi \in C^\infty_0(\overline{\Omega}\,;\R^N)$,
namely they are smooth functions on $\overline{\Omega}$ whose support is compact
in $\overline{\Omega}$.
We assume that $F_i(u,\partial u,\nabla_{\!x}\, \partial u)$ is a smooth function
vanishing to first order at the origin.
Besides, $\partial_0\equiv\partial_t=\partial/\partial t$, 
$\partial_j=\partial/\partial x_j$ ($j=1,2,3$),
$\Delta=\sum_{j=1}^3 \partial_j^2$, $\nabla_{\!x}\, u=(\partial_1 u, \partial_2 u, \partial_3 u)$
and $\partial u=(\partial_t u, \nabla_{\!x}\, u)$. 
In the following we always assume that
\begin{equation}
\frac{\pa F_i}{\pa(\pa_k\pa_\ell u_j)}=\frac{\pa F_j}{\pa(\pa_k\pa_\ell u_i)}
=\frac{\pa F_i}{\pa(\pa_\ell\pa_k u_j)}
\label{symmetric}
\end{equation}
holds for $1\le i, j\le N$ and $1\le k, \ell\le 3$, so that
the hyperbolicity of the system is assured.

First we consider the single speed case (i.e., $c_1=c_2=\cdots =c_N=1$).
If we suppose in addition that quadratic part of the nonlinearity $F_i$ vanishes, then
it was shown in Shibata -- Tsutsumi \cite{ShiTsu86} that 
the mixed problem \eqref{ap1}--\eqref{ap3} admits a unique global small amplitude solution.
Otherwise, in order to get a global existence result, we need a certain
algebraic condition on the nonlinearity in general, due to the blow-up result
for the corresponding Cauchy problem obtained by John \cite{john3} and the finite speed
of propagation. 
One of such conditions is the null condition introduced by Klainerman
\cite{Kla86} (see Definition 1.1 below). 
Under the null condition, Klainerman \cite{Kla86} and 
Christodoulou \cite{Chr86} proved 
global solvability for the Cauchy problem with small initial data
independently by different methods.
This result was extended to the mixed problem by Keel -- Smith -- Sogge \cite{KeSmiSo02G} 
if the obstacle ${\mathcal O}$ is star-shaped, and by Metcalfe \cite{Met04}
if it is non-trapping
(for the case of other space dimensions, we refer to \cite{ShiTsu86}, \cite{Ha95}).

Next we consider the multiple speeds case where the propagation speeds
$c_i$ ($1\le i\le N$) do not necessarily coincide with each other.
Metcalfe -- Sogge \cite{MetSo05} and
Metcalfe -- Nakamura -- Sogge \cite{MetNaSo05a, MetNaSo05b}
extended the global existence result for the mixed problem 
to the multiple speeds case with more general obstacle 
as we shall describe later on
(see \cite{Kov89}, \cite{Yok00}, \cite{Kub-Yok01}, \cite{Kat04:02}, and
 \cite{kayo} for the Cauchy problem in three space dimensions; see also \cite{hk2} for the two space dimensional case).

The aim of this article is to present an alternative approach 
to these works 
which consists of the following two ingredients.
One is the usage of space-time decay estimates for the mixed problem of the linear wave equation
given in Theorem \ref{main3} below,
which directly give us rather detailed decay estimates
\begin{align}\label{std0}
& |u_i(t,x)| \le C \varepsilon (1+t+|x|)^{-1} \log\left(1+\frac{1+c_it+|x|}{1+|c_it-|x|\,|}\right),
\\ \label{std1}
& |\partial u_i(t,x)| \le C \varepsilon (1+|x|)^{-1} (1+|c_i t-|x||)^{-1}
\end{align}
for $(t,x) \in [0,\infty) \times \overline{\Omega}$.
These estimates are refinement of time decay estimates obtained
in the previous works for the mixed problems.
In this way, we do not
need the space--time $L^2$ estimates which has been adopted in the
works \cite{KeSmiSo02G, Met04, MetNaSo05a, MetNaSo05b, MetSo05}.

The other is making use of stronger decay property of a tangential derivative 
to the light cone given in Theorem \ref{D+Es} below.
This idea is recently introduced by the authors \cite{KaKu07},
where the Cauchy problem is studied, and it enables us to deal with 
the null form
without using neither the scaling operator $t\partial_t+x\cdot\nabla_{\!x}\,$ nor Lorentz boost
fields $t \partial_j+x_j \partial_t$\ ($j=1,2,3$).
In this paper, 
we will adopt this approach 
to the mixed problem, and treat the problem without using these vector fields. 
In contrast, the scaling operator has been used in the previous works,
and it makes the argument rather complicated because it does not preserve the Dirichlet boundary condition (\ref{ap2}).
Recently Metcalfe -- Sogge \cite{MetSo07} introduced a simplified approach 
which enables us to use the scaling operator without special care, 
but their approach is applicable only to star-shaped obstacles,
and they assumed that the nonlinearity depends only on derivatives of $u$.



In order to state our result, we need a couple of notions about the obstacle,
the initial data and the nonlinearity.

We remark that we may assume, without loss of generality, that 
${\mathcal O}\subset B_{1}(0)$ by the
scaling and the translation,
where $B_r(z)$ stands for an open ball of
radius $r$ centered at $z \in \R^3$.
Hence we always assume ${\mathcal O}\subset B_1(0)$ in what follows.

Throughout this paper, we denote the standard Lebesgue
and Sobolev spaces by $L^2({\Omega})$ and $H^m({\Omega})$ 
and their norms by $\|\,\cdot : L^2({\Omega})\|$ and 
$\|\,\cdot : H^m({\Omega})\|$, respectively. 
Besides, $H^1_0(\Omega)$ is the completion of
$C^\infty_0({\Omega})$ with respect to $\|\,\cdot : H^1({\Omega})\|$.

\begin{definition}
{\rm (i)}\ 
We say that the obstacle ${\mathcal O}$ is {\bf admissible} if 
there exists a non--negative integer $\ell$ having the following property\,{\rm :}\
Let $v \in C^\infty([0,\infty)\times \overline{\Omega}; \R)$ be a solution
of the homogeneous wave equation 
$(\pa_t^2-c^2\Delta)v=0$ in $[0,\infty)\times \Omega$,
with some constant $c>0$ and the Dirichlet condition,
whose initial value $(v(0,x), (\partial_t v)(0,x))$ vanishes
for $x \in \R^3 \setminus {B_a(0)}$ with some $a>1$.
Then for any $b>1$ we have
\begin{align}\label{obstacle}
&\sum_{|\alpha| \le 1}
 \|\partial^\alpha v(t):L^2({\Omega\,\cap B_b(0)})\|
\\ \nonumber
& \quad \le C \exp(-\sigma t)\,(\|v(0):H^{\ell+1}(\Omega)\|
 {}+\|(\partial_t v)(0):H^{\ell}(\Omega)\|),
\end{align}
where $C$ and $\sigma$ are positive constants depending on $a$, $b$, $c$ and $\Omega$.

\medskip

\noindent
{\rm (ii)}\ 
We say that the initial data $(\phi,\psi)$ satisfies the {\bf compatibility condition}
to infinite order for the mixed problem \eqref{ap1}--\eqref{ap3} 
if the {\rm (}formal{\rm)} 
solution $u$ of the problem satisfies
$(\partial^j_t u)(0,x)=0$ for any
$x \in \partial\Omega$ and 
any non--negative integer $j$
$($notice that the values $(\partial^j_t u)(0,x)$ are determined by $(\phi,\psi)$ and
$F$ successively; for example we have
$\pa_t^2u_i(0,x)=\ve c_i^2 \Delta\phi_i+F_i\bigl(\ve\phi, \ve(\psi, \nabla_x\phi), \ve\nabla_x(\psi, \nabla_x \phi) \bigr)$, and so on$)$.

\medskip

\noindent
{\rm (iii)}\ 
We say that the nonlinearity 
$F=(F_1, F_2, \dots, F_N)$
satisfies the {\bf null condition} associated with the
propagation speeds $(c_1, c_2, \dots, c_N)$ if 
each $F_i$ $(1\le i\le N)$ satisfies
\begin{equation}\label{nullc}
F_i^{(2)}(\lambda, V(\mu,X), W(\nu,X))=0
\end{equation}
for any  
$\lambda$, $\mu$, $\nu \in \Lambda_i$ and 
$X=(X_0, X_1, X_2, X_3)\in \R^{4}$ satisfying 
$X_0^2=c_i^2(X_1^2+X_2^2+X_3^2)$, where $F_i^{(2)}$ is the quadratic part 
of $F_i$, and
$$
 \Lambda_i=\{(\lambda_1, \lambda_2, \ldots, \lambda_N)\in \R^N;
 \lambda_j=0 \text{ if } c_j\ne c_i\}.
$$
Here we put
$V(\mu,X)=(X_a\,\mu_k:\,a=0,1,2,3, \,k=1, \dots, N)$, 
$W(\nu,X)=(X_j X_a \nu_k:\,j=1,2,3, \,a=0,1,2,3, \,k=1, \dots, N)$.
\end{definition}

We often refer to \eqref{obstacle} as the local energy decay. 
We remark that
when ${\mathcal O}$ is non--trapping, the estimate (\ref{obstacle}) holds for $\ell=0$
{\rm (}see for instance Melrose \cite{Mel79}, Shibata -- Tsutsumi \cite{ShiTsu83}{\rm )}.
Even if ${\mathcal O}$ is trapping, it may be admissible in some cases. 
In fact, (\ref{obstacle}) for $\ell=5$ was obtained by Ikawa \cite{Ika82},
provided that ${\mathcal O}$ is a union of disjoint compact sets ${\mathcal O}_1$ and
${\mathcal O}_2$ whose Gaussian curvatures are strictly positive at every point
of their boundaries (see also Ikawa \cite{Ika88}). 

Now we are in a position to state our main result.
\begin{theorem}\label{thm:GE}
Suppose that ${\mathcal O}$ is admissible and that
$(\phi,\psi)$ satisfies the compatibility condition
to infinite order for the problem \eqref{ap1}--\eqref{ap3}.
If $F$ satisfies the null condition associated with $(c_1, c_2, \dots, c_N)$,
then there exists a positive constant $\varepsilon_0$ such that for
all $\varepsilon \in (0,\varepsilon_0)$ the mixed problem \eqref{ap1}--\eqref{ap3}
admits a unique solution $u \in
C^\infty([0,\infty)\times \overline{\Omega}; \R^N)$ satisfying
\eqref{std0} and \eqref{std1}.
\end{theorem}

As we have mentioned in the above, 
the existence part of the Theorem \ref{thm:GE}
 is already known in \cite{MetNaSo05b}
(though the decay property obtained in \cite{MetNaSo05b} is different from ours),
and our aim here is to give a simplified proof for it.

This paper is organized as follows. In the next section we collect notation.
In the section 3 we give some preliminaries needed later on.
The section 4 is devoted to establish pointwise decay estimates. 
Making use of the estimates from the section 4, we give a proof of
Theorem \ref{thm:GE} in the section 5.
\section{Notation}
Let $c>0$.
We shall consider the mixed problem\,:
\begin{align}\label{eq}
& (\partial_t^2-c^2\Delta) v =f, & (t,x) \in (0,T)\times \Omega,
\\ \label{dc}
& v(t,x)=0, & (t,x) \in (0,T)\times \partial\Omega,
\\ \label{id}
& v(0,x)=v_0(x),\ (\partial_t v)(0,x)=v_1(x), & x\in \Omega,
\end{align}
Here $v_0$, $v_1 \in C^\infty_0(\overline{\Omega};\R)$ and 
$f \in C^\infty([0,T)\times \overline{\Omega};\R)$.
We say that $({v}_0, {v}_1, f)$ satisfies the compatibility condition
to infinite order for the problem \eqref{eq}--\eqref{id} if 
$v_j=0$ {on} $\partial\Omega$ for 
any non--negative integer $j$,
where we have set
\begin{equation}\label{data+}
v_j(x)\equiv c^2\Delta v_{j-2}(x)+(\partial_t^{j-2} f)(0,x)
\quad \mbox{for \ $x \in \overline{\Omega}$ \ and \ $j\ge 2$}.
\end{equation}
Let us put $\vec{v}_0:=(v_0, v_1)$ and we denote by $K[\vec{v}_0;c](t,x)$ the solution
of the problem (\ref{eq})--(\ref{id}) with $f \equiv 0$.
While, we denote by $L[f;c](t,x)$ the solution
of the problem (\ref{eq})--(\ref{id}) with $\vec{v}_0 \equiv 0$.

In a similar fashion, putting $\vec{w}_0:=(w_0,w_1)\in C^\infty(\R^3;\R^2)$, we denote by $K_0[\vec{w}_0;c](t,x)$ and
$L_0[g;c](t,x)$ the solution of the following Cauchy problem 
with $g \equiv 0$ and $\vec{w}_0 \equiv 0$, respectively\,:
\begin{align}\label{eq0}
&(\partial_t^2-c^2\Delta) w = g, & (t,x) \in (0,T)\times \R^3,
\\ \label{id0}
& w(0,x)=w_0(x),\ (\partial_t w)(0,x)=w_1(x), & x\in \R^3.
\end{align}

Next we introduce vector fields\,:
\begin{equation}\nonumber 
\partial_0=\partial_t, \quad \partial_j \ (j=1,2,3), \quad
\Omega_{ij}=x_i \partial_j-x_j\partial_i \ (1\le i<j \le 3),
\end{equation}
and we denote them by $Z_j$\,($j=0, 1, \dots, 6$), respectively.
Notice that 
\begin{equation}\label{commute}
[Z_i,\partial_t^2-c^2\Delta]=0 \quad (i=0, 1, \dots, 6),
\end{equation}
where we put $[A,B]:=AB-BA$.
Denoting $Z^\alpha=Z_0^{\alpha_0}Z_1^{\alpha_1}
 \cdots Z_{6}^{\alpha_{6}}$ with a multi--index
$\alpha=(\alpha_0, \alpha_1, \dots, \alpha_{6})$, we set
\begin{equation}\label{norm}
|\varphi (t,x)|_m=\sum_{|\alpha| \le m} |Z^\alpha \varphi(t,x)|,
\quad
\|\varphi(t)\|_m=\|\,|\varphi(t,\cdot)|_m\!:\!{L^2(\Omega)}\|
\end{equation}
for a real or $\R^N$--valued smooth function $\varphi(t,x)$ and 
a non--negative integer $m$.

For $\nu$, $\kappa \in \R$, $c \ge 0$ and  $c_j>0$ ($1\le j\le N$),
we define 
\begin{align}
\label{defPhi}
{\Phi}_\nu(t,x)= &
   \begin{cases}
        \langle t+|x|\rangle^{\nu} & \text{ if } \nu<0,
\\
       \log^{-1}\bigg(2+\displaystyle\frac{\langle t+|x|\rangle}{\langle t-|x|\rangle}\bigg)
         & \text{ if } \nu=0,
\\
       \langle t-|x|\rangle^{\nu}   & \text{ if } \nu>0,  
   \end{cases}
\\
\label{defW}
W_{\nu,\kappa}(t,x)= &
\langle t+|x|\rangle^\nu \Bigl( \min_{0\le j\le N} 
                                \jb{c_jt-|x|} \Bigr)^\kappa,
                                \\
\label{defz}
W^{(c)}_{\nu,\kappa}(t,x)=&
\langle t+|x|\rangle^\nu \Bigl( \min_{0\le j\le N; c_j\ne c} 
          \jb{c_jt-|x|} \Bigr)^\kappa, 
\end{align}
where $c_0=0$ and $\langle y \rangle=\sqrt{1+|y|^2}$ for $y \in \R$ .
We define
\begin{equation}\label{eq:3.5}
 \|g(t)\!:\!{M_k(z)}\| =
   \sup_{(s,x) \in [0,t] \times \R^3} 
     \jb{|x|}\,z(s,x)\,|g(s,x)|_k
\end{equation}
for $t\in [0,T)$, a non--negative integer $k$ and any non--negative
function $z(s,x)$.
Similarly we put
\begin{equation}\label{NfW}
 \|f(t)\!:\!{N_k(z)}\|
=\sup_{(s,x) \in [0,t] \times \Omega} 
     \jb{|x|}\,z(s,x)\,|f(s,x)|_k.
\end{equation}
We also define 
\begin{equation}
B_{\rho, k}[\phi, \psi]=\sup_{y\in \R^3} \jb{|y|}^{\rho} 
\bigl(|\phi(y)|_{k}+|\nabla_x\phi(y)|_k+|\psi(y)|_k\bigr)
\label{HomWei}
\end{equation}
for $\rho\ge 0$, a non--negative integer $k$ and 
$(\phi,\psi) \in (C_0^\infty(\R^3))^2$.

For $a \ge 1$, let $\psi_a$ be a smooth radially symmetric function
on $\R^3$ satisfying
\begin{equation}\label{cutoff}
\psi_a(x)=0 \ (|x| \le a), \quad
\psi_a(x)=1 \ (|x| \ge a+1).
\end{equation}
For $r>0$, we set
$$
\Omega_r=\Omega \cap B_r(0),
$$
where $B_r(x)$ stands for an open ball
of
radius $r$ centered at $x \in \R^3$.

\section{Preliminaries}

First we introduce the local energy decay estimate (\ref{LE})
which works well in getting pointwise estimates for solutions
of our mixed problem.
We also need the elliptic estimate given in Lemma \ref{elliptic}. 
For the completeness, we shall show them in the appendix.

As we have stated in the introduction, we always assume ${\mathcal O}\subset B_1(0)$. 
\begin{lemma}\ \label{local}
Let ${\mathcal O}$ be admissible, and $\ell$ 
be the constant appeared in \eqref{obstacle}.
Suppose that $(\vec{v}_0, f)$ satisfies the compatibility condition 
to infinite order for the mixed problem \eqref{eq}--\eqref{id} and
\begin{eqnarray*}
\text{supp}\,v_j \subset {\Omega_a} \quad (j=0,1), \quad
\text{supp}\,f(t,\cdot) \subset {\Omega_{a}} \quad (t \ge 0)
\end{eqnarray*}
for some $a>1$. 
Let $v$ be the smooth solution of the mixed problem.
Then for any $\gamma>0$, $b>1$ and integer $m$, there exists a
positive constant $C=C(\gamma,a,b,c,m,\Omega)$ such that for $t\in [0,T)$,
\begin{eqnarray}\nonumber
&& \sum_{|\alpha|\le m} \|\partial^\alpha_{t,x} v(t)\!:\!{L^2(\Omega_b)}\|
\le C(1+t)^{-\gamma} \bigg( \|\vec{v}_0\!:\!H^{m+\ell}(\Omega) \times H^{m+\ell-1}(\Omega)\|
\\ \label{LE}
&& \hspace{30mm}  +\sup_{0\le s \le t} (1+s)^\gamma \sum_{|\alpha|
\le m+\ell-1} \|\partial^\alpha_{s,x} f(s)\!:\!{L^2(\Omega)}\|
\bigg).
\end{eqnarray}
\end{lemma}

\vspace{2mm}

\begin{lemma}\label{elliptic}\
Let $\varphi \in H^m(\Omega) \cap H_0^1(\Omega)$ 
for some integer $m(\ge 2)$. 
Then we have
\begin{equation}\label{ap10}
 \| \partial^\alpha \varphi : L^2(\Omega) \| \le
C(\|\Delta \varphi\!:\!{L^2(\Omega)}\| + \|\nabla \varphi\!:\!{L^2(\Omega)}\|)
\end{equation}
for $|\alpha|=m$. 
\end{lemma}

\vspace{2mm}

Next
we 
introduce
a couple of known estimates for
the Cauchy problem. The first one is the decay estimate of
solutions to the homogeneous wave equation, due to Asakura \cite[Proposition 1.1]{asa} 
(observe that the general case can be reduced to the case $m=0$, thanks to (\ref{commute})).
Recall that ${\Phi}_\nu(t,x)$ is the function defined by (\ref{defPhi}).

\begin{lemma}\label{lem:freeH}
Let $c>0$.
For $\vec{w}_0\in (C_0^\infty(\R^3))^2$,
$\rho>0$ and a non--negative integer $m$, there exists a positive constant $C=C(\rho, m, c)$ such that
\begin{equation}\label{decay}
\langle t+|x| \rangle\,{\Phi}_{\rho-1}(ct,x)
   |K_{0}[\vec{w_0}; c](t,x)|_m \le C B_{\rho+1, m}[\vec{w_0}]
\end{equation}
for $(t,x) \in [0,\infty) \times \R^3$.
\end{lemma}

\vspace{2mm}

The second one is the decay estimate for the inhomogeneous wave equation.

\begin{lemma}\label{free}\
Let $c>0$, $\rho>0$, 
and $k$ be a non--negative integer. 
If $\nu=\rho$ and $\kappa>1$, or alternatively if $\nu=\rho+\mu$ and $\kappa=1-\mu$
with some $\mu\in (0,1)$,
then there exists a positive constant $C=C(\nu, \kappa, k, c)$ such that
\begin{equation}\label{ba1}
\langle t+|x| \rangle\,{\Phi}_{\rho-1}(ct,x) |L_{0}[g; c](t,x)|_k
 \le C \|g(t)\!:\!{M_k(W_{\nu,\kappa})}\|
\end{equation}
for $(t,x) \in [0,T) \times \R^3$.
\end{lemma}

\noindent{\it Proof.}\ 
The desired estimate for $k=0$ 
was shown in Theorem 3.4 of Kubota -- Yokoyama \cite{Kub-Yok01}
(see also Lemmas 3.2 and 8.1 in Katayama -- Yokoyama \cite{kayo},
 and Lemma 2.2 in the authors \cite{KaKu07}).

Let $|\alpha| \le k$.
Then it follows from (\ref{commute}) that
\begin{equation}\label{ba11}
Z^\alpha L_0[g; c]=L_0[Z^\alpha g; c]+K_0[(\phi_\alpha,\psi_\alpha); c],
\end{equation}
where we put $\phi_\alpha(x)=(Z^\alpha L_0[g; c])(0,x)$, 
$\psi_\alpha(x)=(\partial_t Z^\alpha L_0[g; c])(0,x)$.
From the equation (\ref{eq0}) we get
$$
\phi_\alpha(x)= \sum_{|\beta| \le |\alpha|-2} C_\beta
(Z^\beta g)(0,x), \quad \psi_\alpha(x)= \sum_{|\beta| \le
|\alpha|-1} C_\beta^\prime (Z^\beta g)(0,x)
$$
with suitable constants $C_\beta$ and $C_\beta^\prime$ 
(cf.~\eqref{data+}). 
Therefore, by virtue of Lemma \ref{lem:freeH}, 
it is enough to show
\begin{equation}\nonumber
\langle t+|x| \rangle\,{\Phi}_{\rho-1}(ct,x) |L_{0}[Z^\alpha g; c](t,x)| 
 \le C \|g(t)\!:\!{M_k(W_{\nu,\kappa})}\|
\end{equation}
for $(t,x) \in [0,T) \times \R^3$.
But this inequality immediately follows from (\ref{ba1}) for $k=0$.
Thus we finish the proof.
\hfill$\qed$

\vspace{2mm}

The third one is the decay estimate of derivatives of solutions to the
inhomogeneous wave equation. 

\begin{lemma}\label{freeD}\
Let $c>0$, and $k$ be a non--negative integer. 

If $\rho=\nu>1$ and $\kappa>1$, or alternatively if
$0<\rho\le 1$, $\nu=1+\mu$ and $\kappa=\rho-\mu$ with some $\mu\in (0, \rho)$,
then there exists a positive constant $C=C(c,\nu,\kappa,k)$
such that
\begin{equation}\label{ba2}
 \langle |x| \rangle \langle ct -|x| \rangle^{\rho}
 |\partial L_{0}[g; c](t,x)|_k
 \le C \|g(t)\!:\!{M_{k+1}(W_{\nu,\kappa})}\|
\end{equation}
for $(t,x) \in [0,T) \times \R^3$.

On the other hand, if $\rho>0$ and $\kappa > 1$, then 
we have
\begin{equation}\label{poc0}
\langle |x| \rangle \langle ct-|x| \rangle^{\rho}
  |\partial L_{0}[g; c](t,x)|_k
 \le C \|g(t)\!:\!{M_{k+1}(W_{\rho,\kappa}^{(c)})}\|
\end{equation}
for $(t,x) \in [0,T) \times \R^3$.
\end{lemma}

\noindent{\it Proof.}\ 
In view of Lemma 3.2 in \cite{Kub-Yok01}, Lemma 8.2 and 
the proof of Lemma 3.2 in \cite{kayo}, we find
that for $0\le a\le 3$,
\begin{align}
\label{KatYok04} 
\jb{|x|}\jb{ct-|x|}^{\rho}|L_0[\pa_a g; c](t,x)|
 \le C \|g(t)\!:\!{M_1(W_{\nu,\kappa})}\|
\end{align}
when  $\rho=\nu>1$ and $\kappa>1$, or when $0<\rho\le1$, $\nu=1+\mu$, and $\kappa=\rho-\mu$ with some $\mu\in (0,\rho)$,
while
\begin{align}
\label{KatYok04-} 
\jb{|x|}\jb{ct-|x|}^{\rho}|L_0[\pa_a g; c](t,x)|
 \le C \|g(t)\!:\!{M_1(W_{\rho,\kappa}^{(c)})}\|,
\end{align}
if $\rho>0$ and $\kappa > 1$
(see also \cite{KaKu07}).

Since $\pa_a L_0[g; c]=L_0[\pa_a g; c]+\delta_{a0}K_0[(0, g(0, \cdot)); c]$ 
for $0\le a\le 3$ with the Kronecker delta $\delta_{ab}$,
\eqref{ba2} and \eqref{poc0} follow from (\ref{ba11}), \eqref{KatYok04}, \eqref{KatYok04-}, and Lemma \ref{lem:freeH}.
This completes the proof.
\hfill$\qed$

\vspace{2mm}

In order to associate these decay estimates with the energy estimate,
we use a variant of the Sobolev type inequality due to Klainerman,
whose proof will be given in the appendix.

\begin{lemma}\label{KlainermanSobolev}\
Let $\varphi \in C_0^2(\overline{\Omega})$.
Then we have
\begin{equation}\label{ap21}
\sup_{x \in \Omega} \jb{|x|} |\varphi(x)|
\le C \sum_{|\alpha| \le 2} \|\widetilde{Z}^\alpha \varphi\!:\!{L^2(\Omega)}\|,
\end{equation}
where $\widetilde{Z}=\{\partial_1,\partial_2,\partial_3,\Omega_{12},\Omega_{23},\Omega_{13} \}$.
\end{lemma}

\vspace{2mm}

Finally, we recall the estimates of the null forms 
from \cite{KaKu07}.
The null forms $Q_0$ and $Q_{ab}$
are defined by
\begin{align}
Q_0(v,w\,;c)=&(\pa_t v)(\pa_t w)-c^2 (\nabla_{\!x}\,v)\cdot (\nabla_{\!x}\,w),\\
Q_{ab}(v,w)=&(\pa_a v)(\pa_b w)-(\pa_b v)(\pa_a w) \quad \text{($0\le a<b\le 3$)}
\end{align}
for a positive constant $c$, and real valued--functions $v=v(t,x)$ and $w=w(t,x)$.
They are closely related to the null condition.

\begin{lemma}\label{NFEs}
Let $c$ be a positive number and $u=(u_1, \dots, u_N)$.
Suppose that $Q$ is one of the null forms.
Then, for a non--negative integer $k$, there exists a positive constant $C=C(c,k)$
such that
\begin{align*}
|Q(u_j, u_k)|_k & \le C
\bigl\{
|\pa u|_{[k/2]} \sum_{|\alpha|\le k} |D_{+,c} Z^\alpha u|
{}+|\pa u|_{k} \sum_{|\alpha|\le [k/2]} |D_{+,c} Z^\alpha u|\\
& \qquad\quad {}+\frac{1}{r}\bigl(|\pa u|_{[k/2]}|u|_{k+1}+|u|_{[k/2]+1}|\pa u|_k
 \bigr)
\bigr\},
\end{align*}
where we put $D_{+,c}=\partial_t+c\,\partial_r$ with $r\pa_r=x\cdot\nabla_{\!x}$ and $r=|x|$. 
\end{lemma}

\section{Basic estimates}

The aim of this section is to establish pointwise decay estimates
for the mixed problem,
which are deduced from corresponding estimates for the Cauchy problem
in combination with the local energy decay.
Theorem \ref{main1} is the result for the homogeneous wave equation, while Theorem \ref{main3} is for the inhomogeneous wave equation.
In order to handle the null forms, we also need some 
estimates, which will be given in Theorem \ref{D+Es}, of a tangential derivative to
the light cone $t=|x|$ which is denoted by $D_{+,c}=\partial_t+c\partial_r$.
To prove these theorems we use
\begin{lemma}\label{KataLem}
Let ${\mathcal O}$ be admissible, and
$\ell$ be the constant in \eqref{obstacle}.
Suppose that $\chi_j$ $(1\le j\le 3)$ are smooth radially symmetric functions
on $\R^3$
satisfying
$$
\supp \chi_1\subseteq B_b(0),\ \supp \chi_2, \supp \chi_3 \subseteq B_a(0),
\ \chi_2=\chi_3 \equiv 0 \text{ on $B_1(0)$}
$$
with some $a(>1)$ and $b(>1)$. 
Let $c>0$, $\nu>0$, $\kappa\ge 0$, and $\kappa_0\ge 0$, 
while $m$ is a non-negative integer.
Then there exists a positive constant $C$ such that
 \begin{align}
 & \langle t \rangle^\nu
   |\chi_1 L[ \chi_2 g;c ](t,x)|_m
\le C \norm{\chi_2 g(t)}{M_{m+\ell+1}(W_{\nu, \kappa})},\label{KataL01}\\
 & \norm{\chi_1L[\chi_2 g; c](t)}{M_m(W_{\nu, \kappa_0})}
      \le C \norm{\chi_2 g(t)}{M_{m+\ell+1}(W_{\nu, \kappa})},
   \label{KataL02}\\
 & \norm{\chi_2L_0[\chi_3 g; c]}{M_m(W_{\nu, \kappa_0})}
\le C\norm{g(t)}{N_m(W_{\nu,\kappa})},
   \label{KataL03} \\
 & \norm{\chi_2K_0[\vec{v}_0; c]}{M_m(W_{\nu, \kappa})}
  \le C B_{\nu+1, m}[\vec{v}_0],
   \label{KataL04} \\
 & \langle t \rangle^\nu |\chi_1K[\chi_2 \vec{v}_0; c](t,x)|_m
  \le C \norm{\vec{v}_0}{H^{m+\ell+2}(\Omega)\times H^{m+\ell+1}(\Omega)},
   \label{KataL05} \\
 & \norm{\chi_1K[\chi_2\vec{v}_0; c](t)}{M_m(W_{\nu, \kappa})}
   \label{KataL06}\\
 & \qquad\qquad\qquad\qquad\qquad   
 \le C \norm{\vec{v}_0}{H^{m+\ell+2}(\Omega) \times H^{m+\ell+1} (\Omega)}
\nonumber 
\end{align}
for any $g\in C^\infty([0,T)\times \Omega)$, and $\vec{v}_0\in  C^\infty_0(\overline{\Omega})$.
\end{lemma}
\noindent{\it Proof.}\ \
First we note that we have
\begin{equation}
|(\chi_1 h)(t,x)|_m\le C \sum_{|\beta|\le m} |\pa_{t, x}^\beta(\chi_1 h)(t,x)|
\label{KataM01}
\end{equation}
for any smooth function $h$ on $[0,T)\times \Omega$,
since $\supp \chi_1\subset B_b(0)$.
We also note that, if $b>0$, $\nu\ge 0$, and $\kappa\ge 0$, then 
$\langle |x|\rangle W_{\nu, \kappa} (t,x)$, 
$\langle{t+|x|}\rangle \Phi_{\nu-1}(ct,x)$, and 
$\langle{t}\rangle^\nu$ are equivalent to each other
for $(t,x)\in [0, \infty)\times B_b(0)$
(observe that we have $W_{\nu, \kappa}(ct,x)\le C\langle t+|x| \rangle^\nu \langle |x|\rangle^\kappa$).

By \eqref{KataM01}, the Sobolev inequality and \eqref{LE} with $\gamma=\nu$, we obtain
\begin{align*}
\langle t \rangle^\nu
|\chi_1 L[ \chi_2 g;c ](t,x)|_m
  \le & C \langle t\rangle^\nu \sum_{|\beta|\le m+2}
        \norm{\pa^\beta L[\chi_2 g;c](t)}{L^2(\Omega_b)}\\
\le & C \sup_{s\in [0,t]} \langle s\rangle^\nu
\sum_{|\beta|\le m+\ell+1} \norm{\pa^\beta (\chi_2g)(s)}{L^2(\Omega)}\\
\le & C \norm{(\chi_2 g)(t)}{M_{m+\ell+1}(W_{\nu,\kappa})},
\end{align*}
which is \eqref{KataL01}.

From \eqref{KataL01}, we find
\begin{align*}
\norm{\chi_1L[\chi_2 g; c](t)}{M_m(W_{\nu, \kappa_0})}\le & 
 C \sup_{(s,x)\in [0,t]\times\R^3}
\langle s\rangle^\nu |\chi_1L[\chi_2g; c](s,x)|_m\\
\le & C\norm{\chi_2 g(t)}{M_{m+\ell+1}(W_{\nu, \kappa})}.
\end{align*}

On the other hand, by \eqref{ba1}, we obtain
\begin{align*}
& \norm{\chi_2L_0[\chi_3g; c](t)}{M_m(W_{\nu, \kappa_0})}\\
& \qquad \le C \sup_{(s,x)\in [0,t]\times \R^3}
 \langle s+|x|\rangle \Phi_{\nu-1}(cs, x) |L_0[\chi_3g; c](s,x)|_m\\
& \qquad \le C \norm{(\chi_3 g)(t)}{M_m(W_{\nu, 2})}
\le C \norm{(\chi_3 g)(t)}{M_m(W_{\nu, \kappa})}.
\end{align*}

Similarly to the proof of \eqref{KataL03}, \eqref{decay}
immediately implies \eqref{KataL04}. 
From \eqref{KataM01}, the Sobolev inequality and \eqref{LE}
we find
\begin{align*}
\langle t\rangle^\nu|\chi_1 K[\chi_2 \vec{v}_0; c](t,x)|_m
\le & C\langle t\rangle^\nu
 \sum_{|\beta|\le m+2}\norm{\pa^\beta K[\chi_2 \vec{v}_0; c](t)}{L^2(\Omega_b)}\\
\le & C\norm{\chi_2 \vec{v}_0}{H^{m+\ell+2}(\Omega)\times H^{m+\ell+1}(\Omega)},
\end{align*}
which leads to \eqref{KataL05}.
Finally, \eqref{KataL06} immediately follows from \eqref{KataL05} in view of
the equivalence of $\langle{|x|}\rangle W_{\nu, \kappa}(t,x)$ 
and $\langle t\rangle^\nu$ in $[0,\infty)\times B_b(0)$.
This completes the proof.\hfill$\qed$
\begin{theorem}\ \label{main1}
Let ${\mathcal O}$ be admissible, 
$\ell$ be the constant in \eqref{obstacle},
and $c>0$.
Suppose that
$\vec{v}_0\in (C_0^\infty(\overline{\Omega}))^2$ and $(\vec{v}_0,0)$ satisfies the
compatibility condition to infinite order for the mixed problem \eqref{eq}--\eqref{id}.
If $\rho>1$ and $k$ is a non--negative integer, then
there exists a constant $C>0$ such that
\begin{equation}\label{m1}
|K[\vec{v}_0; c](t,x)|_k
  \le C\langle t+|x| \rangle^{-1} \langle ct -|x| \rangle^{-(\rho-1)}
B_{\rho+1, k+\ell+3}[\vec{v}_0]
\end{equation}
for $(t,x)\in [0,\infty)\times \Omega$.
\end{theorem}

\noindent{\it Proof.}\ \
First of all, we recall the following representation formula
based on the cut--off method developed by Shibata \cite{Shi83}, and also
by Shibata -- Tsutsumi \cite{ShiTsu86} where $L^p$--$L^q$ time decay estimates for the mixed
problem was obtained (see also \cite{Kub06})\,:
\begin{equation}\label{homo}
K[\vec{v}_0; c](t,x)=\psi_1(x) K_0[\psi_2 \vec{v}_0; c](t,x)
{}+\sum_{i=1}^4 K_i[\vec{v}_0](t,x),
\end{equation}
for $(t,x)\in [0,T)\times \Omega$.
Here $\psi_a$ is defined by (\ref{cutoff}) and we have set
\begin{align}\label{K1}
& K_1[\vec{v}_0](t,x)=(1-\psi_2(x))L\bigl[\,[\psi_1,-c^2\Delta]
K_0[\psi_2 \vec{v}_0;c];c\bigr](t,x),
\\ \label{K2}
&K_2[\vec{v}_0](t,x)\\
\nonumber
& \quad 
=-L_0\bigl[\,[\psi_2,-c^2\Delta]
  L\bigl[\,[\psi_1,-c^2\Delta]
             K_0[\psi_2 \vec{v}_0; c];c\bigr]; c\bigr](t,x),
\\ \label{K3}
& K_3[\vec{v}_0](t,x)=(1-\psi_3(x)) K[(1-\psi_2) \vec{v}_0; c](t,x),
\\ \label{K4}
& K_4[\vec{v}_0](t,x)=-L_0\bigl[\,[\psi_3,-c^2\Delta] K[(1-\psi_2)
\vec{v}_0;c]; c\bigr](t,x).
\end{align}

It is easy to see from (\ref{decay}) for $\rho>1$ that the first term on the
right--hand side of (\ref{homo}) has the desired bound.
Hence our task is to show \eqref{m1} with $K[\vec{v}_0; c]$ replaced by
$K_i[\vec{v}_0]$ ($1\le i\le 4$).

It is easy to check that
\begin{align*}
[\psi_a,-\Delta]u(t,x)= &
   u(t,x) \Delta \psi_a(x)+2\nabla_{\!x}\, u(t,x) \cdot \nabla_{\!x}\, \psi_a(x)\\
=& 2\sum_{j=1}^3 \pa_j\bigl(u(x)\pa_j \psi_a(x)\bigr)-u(x)\Delta\psi_a(x)
\end{align*}
and
\begin{equation*}
\sum_{|\alpha| \le m} \|Z^\alpha
    [\psi_a,-\Delta]u(t)\!:\!L^2(\Omega)\| \le C
\sum_{|\alpha| \le m+1} \|\partial^\alpha u(t)\!:\!L^2(\Omega_{a+1})\|
\end{equation*}
for $t \in [0,T)$, $x \in \Omega$, $a \ge 1$ and any smooth function $u$.
Therefore, by \eqref{KataL01} and \eqref{KataL04} with $\nu=\rho$, we get
\begin{eqnarray*}
|K_1[\vec{v}_0](t,x)|_k \le C\langle t \rangle^{-\rho}
B_{\rho+1, k+\ell+2}[\vec{v_0}],
\end{eqnarray*}
which leads to \eqref{m1} with $K$ replaced by $K_1$,
because
$\text{supp} K_1[\vec{v}_0](t,\cdot) \subset \overline{\Omega_3}$.
On the other hand, \eqref{ba1}, \eqref{KataL02}, and \eqref{KataL04}
with $\nu=\rho$
imply
$$
|K_2[\vec{v}_0](t,x)|_k\le C \langle t+|x|\rangle^{-1} 
\langle ct-|x| \rangle^{-(\rho-1)}
B_{\rho+1, k+\ell+3}[\vec{v}_0].
$$
The bound for $K_3[\vec{v}_0](t,x)$ can be easily obtained by \eqref{KataL05}.
Finally, \eqref{ba1} and \eqref{KataL06} imply the estimate for  
$K_4[\vec{v}_0](t,x)$.
This completes the proof. \hfill$\qed$
\begin{theorem}\ \label{main3}
Let ${\mathcal O}$ be admissible, 
$\ell$ be the constant in \eqref{obstacle},
and $c>0$.
Suppose that
$f \in C^\infty([0,T)\times \Omega)$ and $(0,0,f)$ satisfies the compatibility
condition to infinite order for the mixed problem \eqref{eq}--\eqref{id}.

\noindent 
{\rm (i)}\ Let $\rho>0$. If $\nu=\rho$ and $\kappa>1$, 
or alternatively if $\nu=\rho+\mu$ and $\kappa=1-\mu$
with some $\mu\in (0,1)$,
then there exists a constant $C>0$ such that
\begin{align}\label{ba3}
 \langle t+|x| \rangle {\Phi}_{\rho-1}(ct,x) |L[f; c](t,x)|_k
 \le & C 
 \|f(t)\!:\!{N_{k}(W_{\nu,\kappa})}\|\\
 &+C\norm{f(t)}{N_{k+\ell+3}(W_{\rho, 0})} \nonumber\\
 \le & C 
 \|f(t)\!:\!{N_{k+\ell+3}(W_{\nu,\kappa})}\|
 \nonumber
\end{align}
for $(t,x)\in [0,T)\times \Omega$.
\smallskip\\
{\rm (ii)}\ If $\nu=\rho>1$ and $\kappa>1$, or alternatively
if $0<\rho\le 1$, $\nu=1+\mu$ and $\kappa=\rho-\mu$
with some $\mu\in (0, \rho)$, then we have
\begin{equation}\label{ba4}
 \langle |x| \rangle \langle ct -|x| \rangle^{\rho} |\partial L[f; c](t,x)|_k
 \le C
\|f(t)\!:\!{N_{k+\ell+4}(W_{\nu,\kappa})}\|
\end{equation}
for $(t,x)\in [0,T)\times \Omega$.
\smallskip\\
{\rm (iii)}\ If $\rho>0$ and $\kappa>1$, then we have
\begin{equation}\label{DS}
 \langle |x| \rangle \langle ct -|x| \rangle^{\rho} |\partial L[f; c](t,x)|_k
 \le C
\|f(t)\!:\!{N_{k+\ell+4}(W_{\rho,\kappa}^{(c)})}\|
\end{equation}
for $(t,x)\in [0,T)\times \Omega$.
\end{theorem}

\noindent{\it Proof.}\
Note that $L[f; c]$ has the similar expression to (\ref{homo})\,:
\begin{equation}\label{inhomo}
L[f;c](t,x)=\psi_1(x) L_0[\psi_2 f;c](t,x)+\sum_{i=1}^4 L_i[f](t,x)
\end{equation}
for all $(t,x)\in [0,T)\times \Omega$, where 
\begin{align}\label{L1}
& L_1[f](t,x)=(1-\psi_2(x))L\bigl[\,[\psi_1,-c^2\Delta] L_0[\psi_2 f;c];c
\bigr](t,x),
\\ \label{L2}
& L_2[f](t,x)\\
\nonumber
& \qquad\quad =-L_0\bigl[\,[\psi_2,-c^2\Delta] L\bigl[\,[\psi_1,-c^2\Delta]
L_0[\psi_2 f;c];c \bigr];c \bigr](t,x),
\\ \label{L3}
& L_3[f](t,x)=(1-\psi_3(x)) L[(1-\psi_2) f;c](t,x),
\\ \label{L4}
& L_4[f](t,x)=-L_0\bigl[\,[\psi_3,-c^2\Delta] L[(1-\psi_2) f;c];c\bigr](t,x).
\end{align}
The first term on the right--hand side of (\ref{inhomo}) can be easily treated
by Lemmas \ref{free} and \ref{freeD}.

Let $\rho>0$ and $\kappa\ge 0$.
By \eqref{KataL01} and \eqref{KataL03} with $\nu=\rho$, we obtain
\begin{equation}\label{poD1}
\langle t \rangle^{\rho} |L_i[f](t,x)|_k
 \le C  \|f(t)\!:\!{N_{k+\ell+2}(W_{\rho,\kappa})}\|
\end{equation}
for $i=1,3$. 
It is easy to see that $\langle t+|x|\rangle \Phi_{\rho-1}(ct,x)$ 
and $\langle |x|\rangle \langle ct-|x|\rangle^\rho$ are
equivalent to $\langle t\rangle^\rho$ for $(t,x)\in [0,\infty)\times B_4(0)$.
Therefore, since $\supp L_i[f](t,x)\subset B_4(0)$ for $i=1,3$,
\eqref{poD1} implies the desired estimates
for $L_1[f]$ and $L_3[f]$,
corresponding to \eqref{ba3}, \eqref{ba4} and \eqref{DS}
(note that we also have 
$W_{\rho, \kappa}\le W_{\nu, \kappa}\le W_{\nu, \kappa}^{(c)}$
for $\nu\ge \rho$).

On the other hand, by \eqref{KataL02} and \eqref{KataL03},
we obtain
\begin{align}
\label{KataN01}
\norm{\dal_c L_i[f](t)}{M_{m}(W_{\nu, \kappa_0})}
\le & C \norm{f(t)}{N_{m+\ell+3}(W_{\nu, \kappa})}\\
\le & C \norm{f(t)}{N_{m+\ell+3}(W_{\nu, \kappa}^{(c)})}\ (i=2,4)
\nonumber
\end{align}
for any $\nu>0$, $\kappa_0, \kappa\ge 0$, and $m\ge0$, where $\dal_c=\pa_t^2-c^2\Delta$.
Hence Lemmas \ref{free} and \ref{freeD} imply the desired estimates
for $L_2[f]$ and $L_4[f]$.
This completes the proof. \hfill$\qed$

\vspace{2mm}

\begin{theorem}\label{D+Es}
Let the assumptions in Theorem {\rm \ref{main3}} be fulfilled,
and $1\le \rho\le 2$. 

If $\nu=\rho$ and $\kappa>1$, or
alternatively if $\nu=\rho+\mu$, $\kappa=1-\mu$ with some $\mu\in (0,1)$,
then there exists a positive constant $C=C(\nu,\kappa,c)$ such that
\begin{align}
\label{D+Es-01}
& \jb{|x|} \jb{t+|x|} \jb{ct-|x|}^{\rho-1} 
\sum_{|\alpha|\le k} |D_{+,c}Z^\alpha L[f;c](t,x)|
\\
& \qquad\qquad \le C \log(2+t+|x|)\,\|f(t)\!:\!{N_{k+\ell+5}(W_{\nu,\kappa})}\|.
\nonumber
\end{align}

If $\nu>\rho+1$, we have
\begin{align}
\label{D+Es-02}
& \jb{|x|} \jb{t+|x|} \jb{ct-|x|}^{\rho-1} 
\sum_{|\alpha|\le k} |D_{+,c}Z^\alpha K[\vec{v}_0;c](t,x)|
\\
\nonumber
& \qquad\qquad\qquad\qquad \le C B_{\nu, k+\ell+5}[\vec{v}_0]
\end{align}
for $(t,x)\in [0,T)\times \Omega$.
\end{theorem}
\noindent{\it Proof.}\ 
We consider only (\ref{D+Es-01}), because (\ref{D+Es-02}) can be shown
less hard by using (\ref{m1}).
When $|x| \le 1$, (\ref{D+Es-01}) follows from 
\eqref{ba3} immediately.
While, if $|x|>1$, then we can proceed as in the proof of Theorem 1.2 in \cite{KaKu07},
because ${\mathcal O} \subset B_1(0)$. Here we only give an outline
of the proof.
Setting $U(t, r, \omega)=rL[f; c](t, r\omega)$ for $r>1$ and $\omega\in S^2$,
we have
\begin{align}
 D_{-,c}D_{+,c} U(t,r,\omega)=rf(t,r\omega)
 {}+\frac{c^2}{r}\sum_{1\le j<k\le 3} \Omega_{jk}^2 L[f;c](t, r\omega),
 \label{rad}
\end{align}
where $D_{-,c}=\pa_t-c\pa_r$.
Let $t_0>0$, $r_0>1$ and $\omega_0\in S^2$.
Applying \eqref{ba3} to estimate the second term on 
the right-hand side of \eqref{rad} in terms of 
$\norm{f(t)}{N_{\ell+5}(W_{\nu, \kappa})}$, 
and then integrating the obtained inequality 
along the ray $\{(t, (r_0+c(t_0-t)\omega_0); 0\le t\le t_0\}$
(note that this ray lies in $\Omega$),
we obtain
\begin{align}
\label{Tsuika}
& |D_{+,c}U(t_0, r_0, \omega_0)|\\
& \qquad \le C\jb{t_0+r_0}^{-\rho}\log(2+t_0+r_0)
\norm{f(t_0)}{N_{\ell+5}(W_{\nu, \kappa})}.
\nonumber
\end{align}
Since $rD_{+,c}L[f;c](t, r\omega)=D_{+, c}U(t, r, \omega)-cL[f;c](t, r\omega)$,
\eqref{Tsuika} and \eqref{ba3} imply \eqref{D+Es-01} for $k=0$.
It is easy to obtain \eqref{D+Es-01} for general $k$.
This completes the proof.
\hfill$\qed$

\section{Proof of Theorem 1.2}

In this section we prove Theorem \ref{thm:GE}.
We assume ${\mathcal O}\subset B_{1}(0)$ as before.
Let all the assumptions of Theorem \ref{thm:GE} be fulfilled.

Though there is no essential difficulty in treating the general case
\footnote{
In fact, to treat the general case, we only have to replace the energy inequality 
for the wave equation in Subsections \ref{KEE1}, \ref{KEE2} and \ref{KEE3} 
below with
that for systems of perturbed wave equations
which is also standard (remember that the symmetry conditions \eqref{symmetric}
are assumed).
Such replacement is not needed for pointwise decay estimates,
because loss of derivatives is allowed there.
},
we concentrate on the semilinear case to keep our exposition simple.
Hence we assume $F=F(u, \pa u)$ in what follows.

From the null condition associated with $(c_1, c_2, \ldots, c_N)$,  
we see that
the quadratic part $F_i^{(2)}$ of $F_i$ is independent of $u$, and
can be written as
\begin{equation}
F_i^{(2)}(\pa u)=F_i^{{\rm null}}(\pa u)+R_{I,i}(\pa u)+R_{II,i}(\pa u),
\end{equation}
where
\begin{align*}
F_i^{{\rm null}}(\pa u)=& \sum_{\substack{1\le j,k\le N\\ c_j=c_k=c_i}}
\left(A_{i}^{jk} Q_0(u_j, u_k; c_i)+\sum_{0\le a<b\le 3} B_i^{jk,ab}Q_{ab}(u_j, u_k)
\right),\\
R_{I,i}(\pa u)=&\sum_{\substack{1\le j,k\le N\\ c_j\ne c_k}}\sum_{0\le a, b\le 3}
C_i^{jk,ab} (\pa_a u_j)(\pa_b u_k),\\
R_{II,i}(\pa u)=&\sum_{\substack{1\le j,k\le N\\ c_j=c_k\ne c_i}}\sum_{0\le a, b\le 3}
D_i^{jk,ab} (\pa_a u_j)(\pa_b u_k)
\end{align*}
with suitable constants $A_i^{jk}$, $B_i^{jk, ab}$, $C_i^{jk, ab}$ and $D_i^{jk, ab}$.
We put $$
H_i(u, \pa u)=F_i(u, \pa u)-F_i^{(2)}(\pa u)
$$ for $i=1, 2, \dots, N$, so that $H_i(u,\pa u)=O(|u|^3+|\pa u|^3)$
near $(u, \pa u)=(0,0)$.

Let $u=(u_1, u_2, \dots, u_N)$ be a smooth solution to \eqref{ap1}--\eqref{ap3} on $[0,T)\times \overline{\Omega}$.
We set
\begin{align*}
e_{k,i}[u_i](t,x)=& \jb{t+|x|}\Phi_0(c_it, x)|u_i(t,x)|_{k+1}
{}+\jb{|x|}\jb{c_it-|x|}|\pa u_i(t,x)|_{k}\\
& {}+\frac{\jb{|x|}\jb{t+|x|}}{\log(2+t+|x|)}
\sum_{|\alpha|\le k-1} |D_{+, c_i}Z^\alpha u_i(t,x)|
\end{align*}
for $1\le i \le N$. We also set $e_k[u](t)=\sum_{i=1}^N e_{k,i}[u_i](t,x)$.

We fix $k\ge 6\ell+30$, and
assume that
\begin{equation}
\sup_{0\le t<T} \norm{e_{k}[u](t)}{L^\infty(\Omega)} \le M\ve
\label{InductiveAs}
\end{equation}
holds for some large $M(>1)$ and small $\ve(>0)$, satisfying $M\ve \le 1$.
Since the local existence for the mixed problem has been shown by \cite{ShiTsu86}, what
we need for the proof of the global existence result
is 
a suitable {\it a priori} estimate. 
We will prove that 
\eqref{InductiveAs} implies
\begin{equation}
\label{KFF}
\sup_{0\le t<T} \norm{e_{k}[u](t)}{L^\infty(\Omega)} \le C\ve+CM^2\ve^2.
\end{equation}
From \eqref{KFF} we find that \eqref{InductiveAs} with $M$ replaced by
$M/2$ is true for sufficiently large $M$ and sufficiently small $\ve$,
and the standard continuity argument implies that $e_k[u](t)$ stays bounded
as long as the solution $u$ exists. Theorem \ref{thm:GE}
follows immediately from this {\it a priori} bound.

To this end, the following energy estimate is crucial\,:
\begin{equation}\label{ap16}
 \|\partial u(t)\|_{2k-\ell-8} \le C M\ve (1+t)^{C_* M\ve+\rho_*}
\quad \text{for} \ t \in [0,T),
\end{equation}
where $C$, $C_*$ and $\rho_*$ are positive constants independent of $M$ and $\ve$.
Moreover $\rho_*$ can be chosen arbitrarily small.
In fact, once we find (\ref{ap16}), we can proceed as in the 
case of the corresponding Cauchy problem.
While, unlike the case of the Cauchy problem, 
it is not so simple to get \eqref{ap16},
because of boundary terms coming from the integration--by--parts argument which may cause
some loss of derivatives. 
For this reason, we estimate the space--time gradient
and generalized derivatives separately and improve the estimate of the latter by using the local energy decay.

In the following, we set $r=|x|$. We define
$$
w_-(t,r)=\min_{0\le j\le N} \jb{c_jt-r},\ w_-^{(c)}(t,r)
=\min_{0\le j\le N; c_j\ne c} \jb{c_jt-r}
$$
for $c\ge 0$, with $c_0=0$. 
Note that, for $0\le j,k\le N$, $c_j\ne c_k$ implies
$$
 \jb{c_jt-r}^{-1}\jb{c_kt-r}^{-1}\le C \jb{t+r}^{-1} \min\{\jb{c_jt-r}, \jb{c_kt-r}\}^{-1}.
$$
Notice also that, for any $\mu>0$ and $c>0$, we have
$$
\Phi_0(ct, x)^{-1}\le C \jb{t+r}^\mu \jb{ct-r}^{-\mu},
$$
where $C$ is a positive constant depending only on 
$\mu$ and $c$.

In the arguments below, we always suppose that $M$ is large enough, while
$\ve$ is small enough to satisfy $M\ve<\!\!<1$.
\subsection{Estimates of the energy}\label{KEE1}
First we evaluate the energy involved by time derivatives. 
From \eqref{InductiveAs} we get 
\begin{equation}\nonumber
 |\pa_t^{2k} F^{(2)}(\pa u)(t,x)|
 \le C M\ve \jb{t}^{-1} \sum_{m=0}^{2k} |\pa_t^{m} \pa u(t,x)|,
\end{equation}
and 
\begin{align}\nonumber
& |\pa_t^{2k} H(u, \pa u)(t,x)|\\
& \qquad \le  C|u(t,x)|^3+C
  \,\sum_{m=0}^k \sum_{|\alpha|\le 1}|\pa_t^m \pa_{t,x}^\alpha u(t,x)|^2
      \sum_{m=0}^{2k} |\pa_t^{m} \pa u(t,x)| \nonumber\\
& \qquad \le CM^3\ve^3 \jb{t+r}^{-3+3\mu}w_-(t, r)^{-3\mu} \nonumber\\ \nonumber
& \qquad\quad{}+CM^2\ve^2 \jb{t+r}^{-2+2\mu}w_-(t,r)^{-2\mu}
       \sum_{m=0}^{2k} |\pa_t^{m} \pa u(t,x)|
\end{align}
with small $\mu>0$. Since we have 
$$
 \norm{\jb{t+|\cdot|}^{-3+3\mu}\jb{c_jt-|\cdot|}^{-3\mu}}{L^2(\R^3)}\le 
 C_\mu\jb{t}^{-3/2}
$$
for $\mu>0$ and $0\le j\le N$,
if we set $y(t)= \sum_{m=0}^{2k} \|\pa_t^{m} \pa u(t)\!:\!{L^2(\Omega)}\|$,
then 
we get
$$
\|\pa_t^{2k} F(u, \pa u)(t)\!:\!{L^2(\Omega)}\| 
\le C_0 M\ve (1+t)^{-1} y(t) +CM^3\ve^3(1+t)^{-3/2},
$$
where $C_0$ is a universal constant which is independent of $M$ and $\ve$.
Noting that the boundary condition (\ref{ap2}) implies
$\partial_t^j u(t,x)=0$ for $(t,x)\in [0,T) \times \partial \Omega$
and $0\le j\le 2k+1$, 
we see from the energy inequality 
for the wave equation that 
\begin{equation}\nonumber
\frac{dy}{dt}(t) \le C_0 M\ve (1+t)^{-1} y(t) +CM^3\ve^3(1+t)^{-3/2},
\end{equation}
which yields
\begin{equation}\label{ene1}
\hspace{10mm} y(t) \le (y(0)+CM^3\ve^3) (1+t)^{C_0 M\ve} \le CM\ve (1+t)^{C_0 M\ve}.
\end{equation}

Next we prove that for $0 \le j+m \le 2k$
\begin{equation}\label{ap11}
 \|\partial_t^{j} \nabla_{\!x}\, u(t)\!:\!{H^m(\Omega)}\|
 \le CM\ve (1+t)^{C_0 M\ve}.
\end{equation}
Since (\ref{ap11}) for $m=0$ follows from (\ref{ene1}), it suffices to consider the case $m \ge 1$. 
Then (\ref{ap10}) yields
\begin{equation}\nonumber
\|\pa^\alpha \partial_t^j \nabla_{\!x}\, u(t)\!:\!{L^2(\Omega)}\|
\le C( \|\Delta \partial_t^{j} u(t)\!:\!{{H}^{m-1}(\Omega) }\|
+\|\nabla_{\!x}\, \partial_t^{j} u(t)\!:\!{L^2(\Omega)}\|)
\end{equation}
for $|\alpha|=m$.
Since $0\le j \le 2k-1$, we see from 
\eqref{ap11} for $m=0$ that the second term
is evaluated by $CM\ve (1+t)^{C_0 M\ve}$.
While, using (\ref{ap1}), the first term is estimated by
$$
C( \|\partial_t^{j+2} u(t)\!:\!{{H}^{m-1}(\Omega) }\|
+\|\partial_t^{j} F(u,\pa u)(t)\!:\!{H^{m-1}(\Omega)}\|).
$$
If we set $z_{j,m}(t)= \sum_{s=0}^{j} \|\pa_t^{s} \pa u(t)\!:\!{H^m(\Omega)}\|$,
then we have
$$
\|\pa_t^{j} F(u, \pa u)(t)\!:\!{H^{m-1}(\Omega)}\| 
\le C M\ve (1+t)^{-1} z_{j,m-1}(t) +CM^3\ve^3(1+t)^{-3/2},
$$
as before. 
In conclusion, we get, for $|\alpha|=m$,
\begin{equation}\nonumber
 \|\pa^\alpha \partial_t^{j} \nabla_{\!x}\, u(t)\!:\!{L^2(\Omega)}\|
\le C z_{j+1,m-1}(t) +CM\ve (1+t)^{C_0 M\ve}.
\end{equation}
Since \eqref{ene1} yields $z_{j,0}(t)\le CM\ve (1+t)^{C_0M\ve}$
for $0\le j\le 2k$, we find from the inductive argument in $m(\ge 1)$
that $z_{j,m}(t)\le CM\ve (1+t)^{C_0M\ve}$
for $0\le j+m\le 2k$.
In particular, we obtain (\ref{ap11}).
\subsection{Estimates of the generalized energy, part 1}\label{KEE2}
In this subsection we evaluate the generalized derivatives $\pa Z^\alpha u$ in $L^2(\Omega)$
for $|\alpha| \le 2k-1$.
Fix small $\mu_0>0$.
It follows from (\ref{commute}) that
\begin{eqnarray}\label{ene2}
&& \quad \frac12\frac{d}{dt} 
 \int_{\Omega} \left(|\partial_t Z^\alpha u_i|^2+|\nabla_{\!x}\, Z^\alpha u_i|^2
  \right)\,dx
\\ \nonumber
&& =\int_{\Omega} Z^\alpha F_i(u,\partial u)\,\partial_t Z^\alpha u_i\,dx
 +c_i^2\int_{\partial \Omega} (\nu\cdot \nabla_{\!x}\, Z^\alpha u_i)\,(\partial_t
  Z^\alpha u_i)\,dS,
\end{eqnarray}
where $\nu=\nu(x)$ is the unit outer normal vector at $x \in \partial \Omega$
and $dS$ is the surface measure on $\partial \Omega$.
Observing that $|Z v|\le C \jb{r} |\pa v|$, we obtain
\begin{align}\label{ene3}
\|Z^{\alpha} F(u, \pa u)(t)\!:\!{L^2(\Omega)}\| 
\le & CM\ve (1+t)^{-1} \|\pa u(t)\|_{|\alpha|}\\
&+CM^2\ve^2(1+t)^{-1+2\mu_0}\|\pa u(t)\|_{|\alpha|-1}
\nonumber\\
&+CM^3\ve^3(1+t)^{-3/2}
\nonumber
\end{align}
for $|\alpha| \le 2k-1$ (cf. \eqref{KataKata01} below).

While, since $\partial \Omega \subset B_{1}(0)$, we have
$|Z^\alpha u(t,x)| \le C\sum_{|\beta| \le |\alpha|} |\partial^\beta u(t,x)|$
for $(t,x)\in [0,T) \times \partial \Omega$.
Hence, by the trace theorem, we see that the second term of (\ref{ene2}) is evaluated by
$C \sum_{|\beta| \le |\alpha|+1} \|\partial^\beta \partial u(t)\!:\!{L^2(\Omega_{2})}\|^2$.

Noting that (\ref{ene1}) and (\ref{ap11}) imply  
\begin{equation}
\label{LEK}
 \|\partial^\beta \partial u(t)\!:\!{L^2(\Omega)}\| \le CM\ve (1+t)^{C_0 M\ve}
\text{ for $|\beta| \le 2k$,} 
\end{equation}
we find from (\ref{ene2}) and (\ref{ene3}) that we have
\begin{align*}
\frac{d}{dt}\|\pa u(t)\|_{m}^2 \le &
C_1 M\ve (1+t)^{-1} \|\pa u(t)\|_{m}^2\\
&{}+CM^3\ve^3 (1+t)^{-1+4\mu_0}\|\pa u(t)\|_{m-1}^2+CM^2\ve^2 (1+t)^{2C_0 M\ve}
\end{align*}
for $m\le 2k-1$,
from which we inductively obtain
\begin{equation}
\|\pa u(t)\|_{m}\le CM\ve (1+t)^{C_0M\ve+2\mu_0(m-1)+(1/2)}
\end{equation}
for $m\le 2k-1$,
provided that $\ve$ is so small that $C_1 M\ve \le 1$.
Setting $\gamma=4(k-1)\mu_0$, we obtain
\begin{equation}\label{ap15}
\|\pa u(t)\|_{2k-1}\le CM\ve (1+t)^{C_0M\ve+\gamma+(1/2)}.
\end{equation}
\subsection{Pointwise estimates, part 1}
By (\ref{ap21}) and (\ref{ap15}) we have
\begin{eqnarray}\label{ap21-}
&& \jb{|x|} |\partial u(t,x)|_{2k-3} 
\le C \|\partial u(t)\|_{2k-1} 
 \le  CM\ve (1+t)^{C_0M\ve+\gamma+(1/2)}.
\end{eqnarray}

From \eqref{InductiveAs} we get
\begin{align}
\label{KataKata01}
|F(u, \pa u)(t,x)|_m\le & CM\ve \jb{t+r}^{-1}w_-(t,r)^{-1}
|\pa u(t,x)|_m\\
&{}+CM^2\ve^2\jb{t+r}^{-2+2\mu}w_-(t,r)^{-2\mu}|u(t,x)|_m
\nonumber
\end{align}
for $m\le 2k$ with small $\mu>0$.
We put 
\begin{equation}
\label{KataKata01a}
 U_{m,\lambda}(t)=\sup_{(s, x)\in [0,t]\times \Omega} 
 \sum_{i=1}^N\jb{s+|x|}^{1-\lambda} \Phi_0(c_is, x) |u_i(s,x)|_m
\end{equation}
for $\lambda\ge 0$.
Then \eqref{KataKata01} yields
\begin{align}
\label{KataKata02}
|F(u, \pa u)(t,x)|_m\le & CM\ve \jb{t+r}^{-1}w_-(t,r)^{-1}|\pa u(t,x)|_m\\
&{}+CM^2\ve^2\jb{t+r}^{\lambda-3+3\mu}w_-(t,r)^{-3\mu} U_{m,\lambda}(t)
\nonumber
\end{align}
for $m\le 2k$.
On the other hand,
using $|u(t,x)|_m\le \jb{|x|}|\pa u(t,x)|_{m-1}$ for $m\ge 1$,
and $|u_i(t,x)|\le M\ve \jb{t+r}^{-1+\mu}\jb{c_it-r}^{-\mu}$,
from \eqref{KataKata01} we also obtain
\begin{align}
\label{KataKata03}
|F(u, \pa u)(t,x)|_m\le & CM\ve \jb{t+r}^{-1+2\mu}
w_-(t,r)^{-2\mu}|\pa u(t,x)|_m\\
&{}+CM^3\ve^3\jb{t+r}^{-3+3\mu}w_-(t,r)^{-3\mu}.
\nonumber
\end{align}

Let $\chi$ be a non--negative $C^\infty(\R)$--function satisfying
$\chi(\lambda)=1$ for $\lambda\le 1$, and $\chi(\lambda)=0$ for $\lambda\ge 2$.
We define
\begin{equation}
\label{KataCut01}
 \chi_{c,t_0,x_0}(t,x)=\chi\Bigl(c(t-t_0)+\sqrt{1+|x-x_0|^2}\Bigr)
\end{equation}
for $c>0$ and $(t_0, x_0)\in \Omega$. Then, because of the 
the finite speed of propagation, we have
\begin{equation}
 L[g;c](t_0,x_0)=L[\chi_{c, t_0,x_0}g; c](t_0, x_0).
\label{KataCut02}
\end{equation}
We also have
\begin{equation}
\label{KataCut03}
 \jb{t+|x|}\le C\jb{t_0+|x_0|} 
\end{equation}
for any $(t,x)\in \supp \chi_{c,t_0,x_0}$ with $t\ge 0$, and any 
$(t_0, x_0)\in [0,\infty)\times \Omega$,
where $C$ is a constant depending only on $c$.

Now we set $\lambda=C_0M\ve+2\gamma+(1/2)$.
Using \eqref{ap21-} and \eqref{KataKata02}
with $m=2k-\ell-6$ and $\mu=(1-\gamma)/3$,
we find
\begin{align*}
  &\|\chi_{c_i, t_0, x_0}F_i(u,\partial u)(t_0)\!:\! N_{2k-\ell-6}(W_{1+\gamma, 
1-\gamma})\|\\
   & \qquad \le 
 CM^2\ve^2(1+U_{2k-\ell-6, \lambda}(t_0))
\jb{t_0+|x_0|}^{\lambda}. 
\end{align*}
On the other hand, by \eqref{ap21-} and \eqref{KataKata03} with $m=2k-3$
and $\mu=\gamma/2$, we obtain
  \begin{align*}
   &\|\chi_{c_i, t_0, x_0}F_i(u,\partial u)(t_0)\!:\! N_{2k-3}(W_{1, 0})\|
 \le CM^2\ve^2 \jb{t_0+|x_0|}^{\lambda},
  \end{align*}
since we may assume $2-(3\gamma/2)\ge 1$.

In view of \eqref{KataCut03}, by using \eqref{m1} and the first inequality in \eqref{ba3} with $(\rho, \nu, \kappa)=(1, 1+\gamma, 1-\gamma)$, we obtain
\begin{align*}
& U_{2k-\ell-6, \lambda}(t)\le C\ve+CM^2\ve^2
(1+U_{2k-\ell-6, \lambda}(t))
\end{align*}
with $\lambda=C_0M\ve+2\gamma+(1/2)$,
which leads to
\begin{equation}
 \sum_{i=1}^N \jb{t+|x|}^{(1/2)-C_0M\ve-2\gamma}\Phi_0(c_it, x)
 |u_i(t,x)|_{2k-\ell-6}\le CM\ve 
\label{KataKata04}
\end{equation}
for $(t,x)\in [0, T)\times \Omega$, since we may assume $CM^2\ve^2\le 1/2$.
\subsection{Estimates of the generalized energy, part 2}\label{KEE3}
Since $\Phi_0(c_it,x)$ is bounded for $(t, x)\in [0,\infty)\times \Omega_2$,
from \eqref{KataKata04} we get
\begin{align}
\label{ap17}
\norm{|u(t)|_{2k-\ell-6}}{L^2(\Omega_2)}
\le & C\norm{|u(t)|_{2k-\ell-6}}{L^\infty(\Omega_2)}\\
\le & CM\ve \jb{t}^{-(1/2)+C_0M\ve+2\gamma}, \nonumber
\end{align}
instead of \eqref{LEK}.
Now (\ref{ene2}), (\ref{ene3}) and (\ref{ap17}) yield
\begin{align}\nonumber
\frac{d}{dt}\|\pa u(t)\|_{m}^2 \le &
C_2 M\ve (1+t)^{-1} \|\pa u(t)\|_{m}^2 
\\ \nonumber
&+CM^3\ve^3(1+t)^{-1+4\mu_0} \|\pa u(t)\|_{m-1}^2
\nonumber\\
&+CM^2\ve^2 (1+t)^{-1+4\gamma+2C_0 M\ve},
\nonumber
\end{align}
for $m\le 2k-\ell-8$,
which inductively leads to (\ref{ap16}) with $C_*=C_0+C_2/2$ and $\rho_*=4\gamma$.
\subsection{Pointwise estimates, part 2}
(\ref{ap21}) and (\ref{ap16}) imply
\begin{align}\label{KataKata05}
& \jb{|x|} |\partial u(t,x)|_{2k-\ell-10} 
\le  CM\ve (1+t)^{\delta}
\end{align}
for $0<\ve<\rho_*/(C_*M)$, where we have set $\delta=2\rho_*$.
Note that we can take $\rho_*$ arbitrarily small, hence we may assume that
$\delta$ is small enough in the following.

Using \eqref{KataKata05} and \eqref{KataKata02}
with $m=2k-2\ell-13$, and $\mu=(1-\delta)/3$,
we find
\begin{align*}
  &\|\chi_{c_i, t_0, x_0}F_i(u,\partial u)(t_0)\!:\! N_{2k-2\ell-13}
 (W_{1+\delta, 1-\delta})\|\\
   & \qquad \le 
 CM^2\ve^2(1+U_{2k-2\ell-13, 2\delta}(t_0))
\jb{t_0+|x_0|}^{2\delta}. 
\end{align*}
On the other hand, by \eqref{KataKata05} and \eqref{KataKata03} with $m=2k-\ell-10$
and $\mu=\delta/3$, we obtain
  \begin{align*}
   &\|\chi_{c_i, t_0, x_0}F_i(u,\partial u)(t_0)\!:\! N_{2k-\ell-10}(W_{1, 0})\|
 \le CM^2\ve^2 \jb{t_0+|x_0|}^{2\delta},
  \end{align*}
since we may assume $2-\delta\ge 1$.
Now, similarly to \eqref{KataKata04}, these estimates end up with
\begin{equation}
 \sum_{i=1}^N \jb{t+|x|}^{1-2\delta}\Phi_0(c_it, x)
 |u_i(t,x)|_{2k-2\ell-13}\le CM\ve 
\label{KataKata06}
\end{equation}
for $(t,x)\in [0, T)\times \Omega$.

From \eqref{KataKata02} (with $\mu=(1+\delta)/3$), \eqref{KataKata05} and \eqref{KataKata06},
we get
\begin{align}
&
\|\chi_{c_i, t_0, x_0}F_i(u,\partial u)(t_0)\!:\! N_{2k-2\ell-13}(W_{1+\delta, 1+\delta})\|
\\ 
\nonumber & \qquad\qquad\qquad\qquad\qquad
\le CM^2\ve^2 \jb{t_0+|x_0|}^{4\delta}.
\end{align}

From \eqref{m1}, \eqref{ba4}, \eqref{D+Es-01} and \eqref{D+Es-02},
we obtain
\begin{align}
& \jb{r}\jb{t+r}^{-4\delta}\jb{c_it-r}^{1+\delta}
|\pa u_i(t,x)|_{2k-3\ell-17}\le CM\ve,
\label{KataKata07}\\
& \jb{r}\jb{t+r}^{1-5\delta}\jb{c_it-r}^\delta
\sum_{|\alpha|\le 2k-3\ell-18}|D_{+, c_i} Z^\alpha u_i(t,x)|
\le CM\ve
\label{KataKata08}
\end{align}
for $1\le i\le N$ and $(t,x)\in [0,T)\times \Omega$,
where we have used $\log(2+t+r)\le C\jb{t+r}^{\delta}$.
\subsection{Pointwise estimates, part 3}
From now on, we take advantage of detailed structure of our nonlinearity.

Note that 
$r$ is equivalent to $\jb{t+r}$,
when $r\ge 1$ and $|c_it-r|<(c_it/2)$.
By Lemma \ref{NFEs}, with the help of
\eqref{InductiveAs}, \eqref{KataKata06}, \eqref{KataKata07},
and \eqref{KataKata08}, we obtain
\begin{align}
|F_i^{{\rm null}}(\pa u)(t,x)|_{2k-3\ell-18} 
\le CM^2\ve^2 \jb{t+r}^{-3+5\delta}
\jb{c_it-r}^{-1-\delta}
\end{align}
for $(t,x)$ satisfying $r\ge 1$ and $|c_it-r|<(c_it/2)$.

On the other hand, $\jb{c_it-r}$ is equivalent to $\jb{t+r}$,
when $r<1$ or $|c_it-r|\ge (c_it/2)$. Hence, observing that $F_i^{{\rm null}}$
is quadratic with respect to $\pa u$, 
from \eqref{InductiveAs} and \eqref{KataKata07} we get
\begin{equation}
|F_i^{{\rm null}}(\pa u)(t,x)|_{2k-3\ell-18}
\le CM^2\ve^2 \jb{t+r}^{-2+3\delta}
\jb{r}^{-2}
\end{equation}
for $(t,x)$ satisfying $r<1$ or $|c_it-r|\ge (c_it/2)$.

Now we find
\begin{align}
\label{KataKata11}
\norm{F_i^{{\rm null}}(\pa u)(t)}{N_{2k-3\ell-18}(W_{\nu, \kappa})}
\le CM^2\ve^2
\end{align}
with some $\nu>1$ and $\kappa>1$, since we may assume
$2-5\delta>1$.

\eqref{InductiveAs} and \eqref{KataKata07} yield
\begin{align}
& |R_{I,i}(\pa u)(t,x)|_{2k-3\ell-18}\\
& \qquad \le CM^2\ve^2 \jb{r}^{-2}\jb{t+r}^{4\delta}
\sum_{c_j\ne c_k}\jb{c_jt-r}^{-1}\jb{c_kt-r}^{-1-\delta}
\nonumber\\
& \qquad \le CM^2\ve^2\jb{r}^{-1}\jb{t+r}^{-2+4\delta}
w_-(t,r)^{-1-\delta}
\nonumber
\end{align}
for $(t,x)\in [0, T)\times \Omega$ with $c_0=0$.
Since we may assume $2-4\delta>1$, we obtain
\begin{align}
\label{KataKata13}
\norm{R_{I,i}(\pa u)(t)}{N_{2k-3\ell-18}(W_{\nu, \kappa})}\le CM^2\ve^2
\end{align}
with some $\nu>1$ and $\kappa>1$.

Similarly, we have
\begin{align}
|R_{II,i}(\pa u)(t,x)|_{2k-3\ell-18}
& \le CM^2\ve^2\jb{r}^{-1}\jb{t+r}^{-1+4\delta}
\\
& \qquad\qquad\qquad \times
w_-^{(c_i)}(t,r)^{-2-\delta},
\nonumber
\end{align}
which yields
\begin{align}
\label{KataKata15}
\norm{R_{II,i}(\pa u)(t)}{N_{2k-3\ell-18}(W_{-1+4\delta, \kappa}^{(c_i)})}\le CM^2\ve^2
\end{align}
with some $\kappa>1$.

From \eqref{InductiveAs}, \eqref{KataKata06} and \eqref{KataKata07}
we have
\begin{align}
\label{KataKata16}
& |H_i(u, \pa u)(t,x)|_{2k-3\ell-18}\\
& \qquad \le CM^3\ve^3
\jb{t+r}^{-3+3\mu+4\delta}
w_-(t,r)^{-3\mu}
\nonumber
\end{align}
with small $\mu>0$,
which implies
\begin{equation}
\label{KataKata17}
\norm{H_i(u,\pa u)(t)}{N_{2k-3\ell-18}(W_{1+\delta, (1-4\delta)-\delta})}
\le CM^2\ve^2.
\end{equation}

Finally,
\eqref{ba3}, \eqref{ba4} and \eqref{D+Es-01} lead to
\begin{equation}
e_{2k-4\ell-22, i}\bigl[L[F_i^{{\rm null}}+R_{I,i}; c_i]\bigr](t,x)
\le CM^2\ve^2
\label{KataF01}
\end{equation}
in view of \eqref{KataKata11} and \eqref{KataKata13}.
On the other hand, \eqref{KataKata15} and \eqref{DS} yield
\begin{equation}
\label{KataF02}
\jb{r}\jb{c_it-r}^{1-4\delta}
|\pa L[R_{II,i}; c_i](t,x)|_{2k-4\ell-22}\le CM^2\ve^2,
\end{equation}
while \eqref{KataKata17} and \eqref{ba4} with $(\rho,\nu, \kappa)=(1-4\delta,
1+\delta, (1-4\delta)-\delta)$ imply
\begin{equation}
\label{KataF03}
\jb{r}\jb{c_it-r}^{1-4\delta}
|\pa L[H_{i}; c_i](t,x)|_{2k-4\ell-22}\le CM^2\ve^2.
\end{equation}

From \eqref{KataF01}, \eqref{KataF02} and \eqref{KataF03}, 
together with \eqref{m1}, we obtain
\begin{equation}
\jb{r}\jb{c_it-r}^{1-4\delta}|\pa u_i(t,x)|_{2k-4\ell-22}\le 
CM\ve.
\label{KataF04}
\end{equation}
\subsection{Pointwise estimates, the final part}
By \eqref{InductiveAs} and \eqref{KataF04}, we obtain
\begin{align}
& |R_{II,i}(\pa u)(t,x)|_{2k-4\ell-22}\\
& \qquad \le CM^2\ve^2\jb{r}^{-1}\jb{t+r}^{-1}
w_-^{(c_i)}(t,r)^{-2+4\delta},
\nonumber
\end{align}
which leads to
\begin{align}
\norm{R_{II,i}(\pa u)(t)}{N_{2k-4\ell-22}(
W_{1, \kappa}^{(c_i)})}\le CM^2\ve^2
\end{align}
with some $\kappa>1$, since we may assume $2-4\delta>1$.
Hence \eqref{ba3}, \eqref{DS} and \eqref{D+Es-01}
imply
\begin{equation}
e_{2k-5\ell-26,i}\bigl[ L[R_{II,i}; c_i] \bigr]
(t,x)\le CM^2\ve^2
\label{KataF06}
\end{equation}
(observe that we have $W_{1,\kappa}\le W_{1, \kappa}^{(c_i)}$).

By \eqref{InductiveAs} and \eqref{KataF04}, we also obtain
\begin{align}
& |H_i(u,\pa u)(t,x)|_{2k-5\ell-26}\\
& \quad \le 
CM^3\ve^3 \jb{r}^{-1} \jb{t+r}^{-2+2\mu}
w_-(t,r)^{-1+4\delta-2\mu}
\nonumber\\
&\quad \qquad {}+CM^2\ve^2\jb{t+r}^{-3+3\mu}
w_-(t,r)^{-3\mu}
U_{2k-5\ell-26, 0}(t)
\nonumber
\end{align}
with small $\mu>0$,
where $U_{m, \lambda}$ is given by \eqref{KataKata01a}.
Since we may assume $-1+4\delta<0$, we have
\begin{align}
\label{KataF09}
& \norm{H_i(u, \pa u)(t)}{N_{2k-5\ell-26}(W_{1+\mu, 1-\mu})} \\
& \qquad \le CM^2\ve^2(M\ve+U_{2k-5\ell-26,0}(t))
\nonumber
\end{align}
From \eqref{KataKata16} we also have
\begin{equation}
\norm{H_i(u, \pa u)(t)}{N_{2k-4\ell-23}(W_{1, 0})} 
\le CM^3\ve^3.
\end{equation}
Now the first inequality in \eqref{ba3} leads to
\begin{align}
\label{KataF07} 
& \jb{t+r}\Phi_0(c_it, x)|L[H_i; c_i](t,x)|_{2k-5\ell-26}\\
& \qquad\qquad\qquad \le
CM^2\ve^2(M\ve+U_{2k-5\ell-26,0}(t)).
\nonumber
\end{align}
\eqref{KataF01}, \eqref{KataF06} and \eqref{KataF07} imply
$$
U_{2k-5\ell-26,0}(t)\le C\ve+CM^2\ve^2(1+U_{2k-5\ell-26,0}),
$$
which yields
\begin{equation}
\label{KataF08}
 \jb{t+r}\Phi_0(c_it, x)|u_i(t,x)|_{2k-5\ell-26}
 \le C\ve+CM^2\ve^2,
\end{equation}
provided that $\ve$ is sufficiently small.
In view of \eqref{KataF09} and \eqref{KataF08},
we obtain
$$
 \norm{H_i(u, \pa u)(t)}{N_{2k-5\ell-26}(W_{1+\mu, 1-\mu})} 
 \le CM^3\ve^3.
$$
Now \eqref{ba4} and \eqref{D+Es-01} with $(\rho, \nu, \kappa)=(1, 1+\mu, 1-\mu)$ 
imply
\begin{align}
\label{KataF10}
& \jb{r}\jb{c_it-r}|\pa L[H_i; c_i](t,x)|_{2k-6\ell-30}\le CM^3\ve^3,\\
\label{KataF11}
& \frac{\jb{r}\jb{t+r}}{\log(2+t+r)}\sum_{|\alpha|\le 2k-6\ell-31}
|D_{+, c_i}Z^\alpha L[H_i; c_i](t,x)|\le CM^3\ve^3.
\end{align}

Finally, since $2k-6\ell-30\ge k$,
from \eqref{KataF01}, \eqref{KataF06}, \eqref{KataF08},
\eqref{KataF10} and \eqref{KataF11}, we obtain 
\eqref{KFF}.
This completes the proof.
\hfill$\qed$
\subsection{Concluding remark}
If we consider the single speed case
$c_1=c_2=\cdots c_N=1$, we can replace $e_{k}[u](t)$ by
\begin{align*}
\widetilde{e}_k[u](t,x)=&\jb{t+|x|}\jb{t-|x|}^\rho|u(t,x)|_{k+1}
{}+\jb{|x|}\jb{t-|x|}^{1+\rho}|\pa u(t,x)|_{k}\\
& {}+\frac{\jb{|x|}\jb{t+|x|}\jb{t-|x|}^\rho}{\log(2+t+|x|)}
\sum_{|\alpha|\le k-1} |D_{+, 1}Z^\alpha u(t,x)|
\end{align*}
with some $\rho \in (1/2, 1)$ as in the Cauchy problem treated in \cite{KaKu07},
and we can show $\norm{\widetilde{e}_k[u](t)}{L^\infty(\R^3)}\le M\ve$
for $0\le t<\infty$. The proof becomes much simpler because of the better
decay of the solution.
\renewcommand{\theequation}{A.\arabic{equation}}

\setcounter{equation}{0}  
\renewcommand{\thelemma}{A.\arabic{lemma}}

\renewcommand{\thetheorem}{A.\arabic{theorem}}

\setcounter{theorem}{0}
\section*{Appendix}

\noindent{\it Proof of Lemma \ref{elliptic}.}\ 
We shall show (\ref{ap10}) only for $m=2$, because the general case can be obtained analogously by the inductive argument.
Let $\chi$ be a $C^\infty_0(\R^3)$ function such that
$\chi \equiv 1$ in a neighborhood of ${\mathcal O}$.
Let $\text{supp}\,\chi \subset B_R(0)$ for some $R>1$.
We set $\varphi_1=\chi \varphi$ and $\varphi_2=(1-\chi) \varphi$, so that $\varphi=\varphi_1+\varphi_2$.

First we prove, for $|\alpha|=2$,
\begin{equation}\label{ap10bis}
 \| \pa^\alpha \varphi_2\!:\!{L^2(\Omega)}\| \le
C(\|\Delta \varphi\!:\!{L^2(\Omega)}\| 
  +\|\nabla \varphi\!:\!{L^2(\Omega)}\|).
\end{equation}
Since $\| \pa^\alpha w\!:\!{L^2(\R^3)}\| \le
C\|\Delta w\!:\!{L^2(\R^3)}\|$ for $|\alpha|=2$ and 
$w \in H^2(\R^3)$, the left--hand side of \eqref{ap10bis} is estimated by
\begin{equation}\nonumber
 C\|\Delta \varphi_2\!:\!{L^2(\Omega)}\|
\le C(\|\varphi\!:\!{L^2(\Omega_R)}\|
  +\|\nabla \varphi\!:\!{L^2(\Omega)}\|
  +\|\Delta \varphi\!:\!{L^2(\Omega)}\|).
\end{equation}
Thanks to the estimate
\begin{equation}\label{ha}
  \|w\!:\!{L^2(\Omega_R)}\| \le C R^2
  \|\nabla w\!:\!{L^2(\Omega)}\|
\end{equation}
for $w \in H_0^1(\Omega)$
(for the proof, see \cite{LaPh}), we obtain \eqref{ap10bis}.

Next we estimate $\varphi_1$.
We shall use the following well--known elliptic estimate
(see Chapter 9 in \cite{GiTr} for instance):
\begin{eqnarray}\nonumber
 \|w\!:\!{H^{k+2}(\Omega_R)}\| \le
C(\|\Delta w\!:\!{H^k(\Omega_R)}\| 
  +\|w\!:\!{L^2(\Omega_R)}\|)
\end{eqnarray}
for $w \in H^{k+2}(\Omega_R) \cap H^1_0(\Omega_R)$ with a non--negative integer $k$.

Since $\text{supp}\,\chi \subset B_R(0)$, we have $\varphi_1 \in H_0^1(\Omega_R)$.
Therefore, the application of the above estimate for $k=0$ in 
combination with (\ref{ha}) gives
\begin{equation}\label{ap10bis2}
 \| \varphi_1\!:\!{H^2(\Omega)}\| \le
C(\|\Delta \varphi\!:\!{L^2(\Omega)}\| 
  +\|\nabla \varphi\!:\!{L^2(\Omega)}\|).
\end{equation}
Thus (\ref{ap10}) for $m=2$ follows from (\ref{ap10bis}) and (\ref{ap10bis2}). 
\hfill$\qed$

\vspace{2mm}

\noindent{\it Proof of Lemma \ref{local}.}\ 
If $v$ is the smooth solution of the mixed problem \eqref{eq}--\eqref{id}, then
it follows that 
\begin{equation}\nonumber 
\partial_t^j v(t,x)=K[(v_j,v_{j+1}); c](t,x)+
\int_0^t K[(0,\partial_s^j f(s)); c](t-s,x) ds
\end{equation}
for any non--negative integer $j$ and any $(t,x) \in [0,T) \times \Omega$, 
where $v_j$ are given by (\ref{data+}).
By (\ref{obstacle}) we have, for ${|\alpha| \le 1}$,
\begin{eqnarray}\label{obstacle1}
&& \hspace{4mm} 
 \|\partial^\alpha K[(v_j,v_{j+1}); c](t):L^2({\Omega_b})\|
\\ \nonumber
&& \le C \exp(-\sigma t)\,(\|v_j:H^{\ell+1}(\Omega)\|+\|v_{j+1}:H^{\ell}(\Omega)\|)
\\ \nonumber
&& \le C \exp(-\sigma t)\,(\|v_0:H^{\ell+j+1}(\Omega)\|+\|v_{1}:H^{\ell+j}(\Omega)\|
\\ \nonumber
&& \hspace{40mm}
   +\sum_{|\alpha| \le \ell+j-1} \| (\partial_{s,x}^\alpha f)(0) : L^2(\Omega)\|)
\end{eqnarray}
and
\begin{eqnarray}\label{obstacle2}
&& \hspace{4mm} 
 \int_0^t \|\partial^\alpha K[(0,\partial_s^j f(s)); c](t-s) :L^2({\Omega_b})\| ds
\\ \nonumber
&& \le C \int_0^t \exp(-\sigma (t-s))\,\| \partial_s^j f(s) :H^{\ell}(\Omega)\| ds
\\ \nonumber
&& \le C(1+t)^{-\gamma} \sup_{0\le s \le t} (1+s)^\gamma 
  \| \partial_s^j f(s) :H^{\ell}(\Omega)\| 
\end{eqnarray}
for any $\gamma>0$.
Therefore for ${|\alpha| \le 1}$ and any non--negative integer $j$, 
we have
\begin{eqnarray}\label{LE1}
&& 
 \| \partial^\alpha \partial^j_{t} v(t)\!:\!{L^2(\Omega_b)}\|
\le C(1+t)^{-\gamma} \,( \|\vec{v}_0\!:\!H^{\ell+j+1}(\Omega) \times H^{\ell+j}(\Omega)\|
\\ \nonumber
&&  \hspace{40mm}
 +\sum_{|\alpha| \le \ell+j} \sup_{0\le s \le t} (1+s)^\gamma 
\|\partial^\alpha_{s,x} f(s)\!:\!{L^2(\Omega)}\|).
\end{eqnarray}

In order to evaluate $\partial^\alpha v$ for ${|\alpha| \le m}$, 
we have only to combine (\ref{LE1}) with a variant of (\ref{elliptic})\,:
\begin{equation}\label{LE2}
 \|\varphi\!:\!{H^m(\Omega_b)}\| \le
C(\|\Delta \varphi\!:\!{H^{m-2}(\Omega_{b^\prime})}\| +\|\varphi\!:\!{H^1(\Omega_{b^\prime})}\|),
\end{equation}
where $1<b<b^\prime$ and $\varphi \in H^m(\Omega) \cap 
H_0^1(\Omega)$ with $m \ge 2$. 
This completes the proof.
\hfill$\qed$

\vspace{2mm}

\noindent{\it Proof of Lemma \ref{KlainermanSobolev}.}\ It is well-known that for $w \in
C_0^2(\R^3)$ we have
\begin{eqnarray}\nonumber
\sup_{x \in \R^3} |x||w(x)| \le C \sum_{|\alpha| \le 2}
\|\widetilde{Z}^\alpha w\!:\!{L^2(\R^3)}\|
\end{eqnarray}
(for the proof, see e.g. \cite{kl0}). Rewriting $\varphi$ as
$\varphi=\psi_1 \varphi+(1-\psi_1) \varphi$ with $\psi_1$ in (\ref{cutoff}), we see that the left--hand side on
(\ref{ap21}) is evaluated by
\begin{eqnarray}\nonumber
&& \hspace{4mm}
 C \sup_{x \in \R^3} |x| |\psi_1(x) \varphi(x)|
+C \sup_{x \in \Omega} |(1-\psi_1(x)) \varphi(x)|
\\ \nonumber
&&  \le
 C \sum_{|\alpha| \le 2} \|\widetilde{Z}^\alpha (\psi_1 \varphi)\!:\!{L^2(\R^3)}\|
+C \sum_{|\alpha| \le 2} \|\partial^\alpha ((1-\psi_1) \varphi)\!:\!{L^2(\Omega_2)}\|
\\ \nonumber
&&  \le C \sum_{|\alpha| \le 2} \|\widetilde{Z}^\alpha \varphi\!:\!{L^2(\Omega)}\|,
\end{eqnarray}
hence we obtain (\ref{ap21}). This completes the proof.
\hfill$\qed$
\medskip\\
\begin{center}
{\sc Acknowledgments}
\end{center}
The authors would like to express their gratitude to
Prof.~S.~Alinhac for his useful comments on 
the preliminary version of this paper.

\end{document}